\documentclass[oneside,english]{amsart}
\usepackage[T1]{fontenc}
\usepackage[latin9]{inputenc}
\usepackage{amstext}
\usepackage{amsthm}
\usepackage{amssymb}
\usepackage{graphicx}

\makeatletter
\numberwithin{equation}{section}
\numberwithin{figure}{section}

\date{Created December 2019; updated \today.}

\usepackage{cite}
\usepackage{url}

\makeatother

\usepackage{babel}
\begin{document}
\title[Solar-weighted Gaussian quadrature]{Accurate solar-power integration: Solar-weighted Gaussian quadrature}
\author{Steven G. Johnson, MIT Applied Mathematics}
\address{Department of Mathematics, Massachusetts Institute of Technology,
Cambridge~MA~02139.}
\begin{abstract}
In this technical note, we explain how to construct Gaussian quadrature
rules for efficiently and accurately computing integrals of the form
$\int S(\lambda)f(\lambda)d\lambda$ where $S(\lambda)$ is the solar
irradiance function tabulated in the ASTM standard and $f(\lambda)$
is an arbitary application-specific smooth function. This allows the
integral to be computed accurately with a relatively small number
of $f(\lambda)$ evaluations despite the fact that $S(\lambda)$ is
non-smooth and wildly oscillatory. Julia software is provided to compute
solar-weighted quadrature rules for an arbitrary bandwidth or number
of points. We expect that this technique will be useful in solar-energy
calculations, where $f(\lambda)$ is often a computationally expensive
function such as an absorbance calculated by solving Maxwell's equations.
\end{abstract}

\maketitle

\section{Introduction}

\begin{figure}
\begin{centering}
\includegraphics[width=0.7\columnwidth]{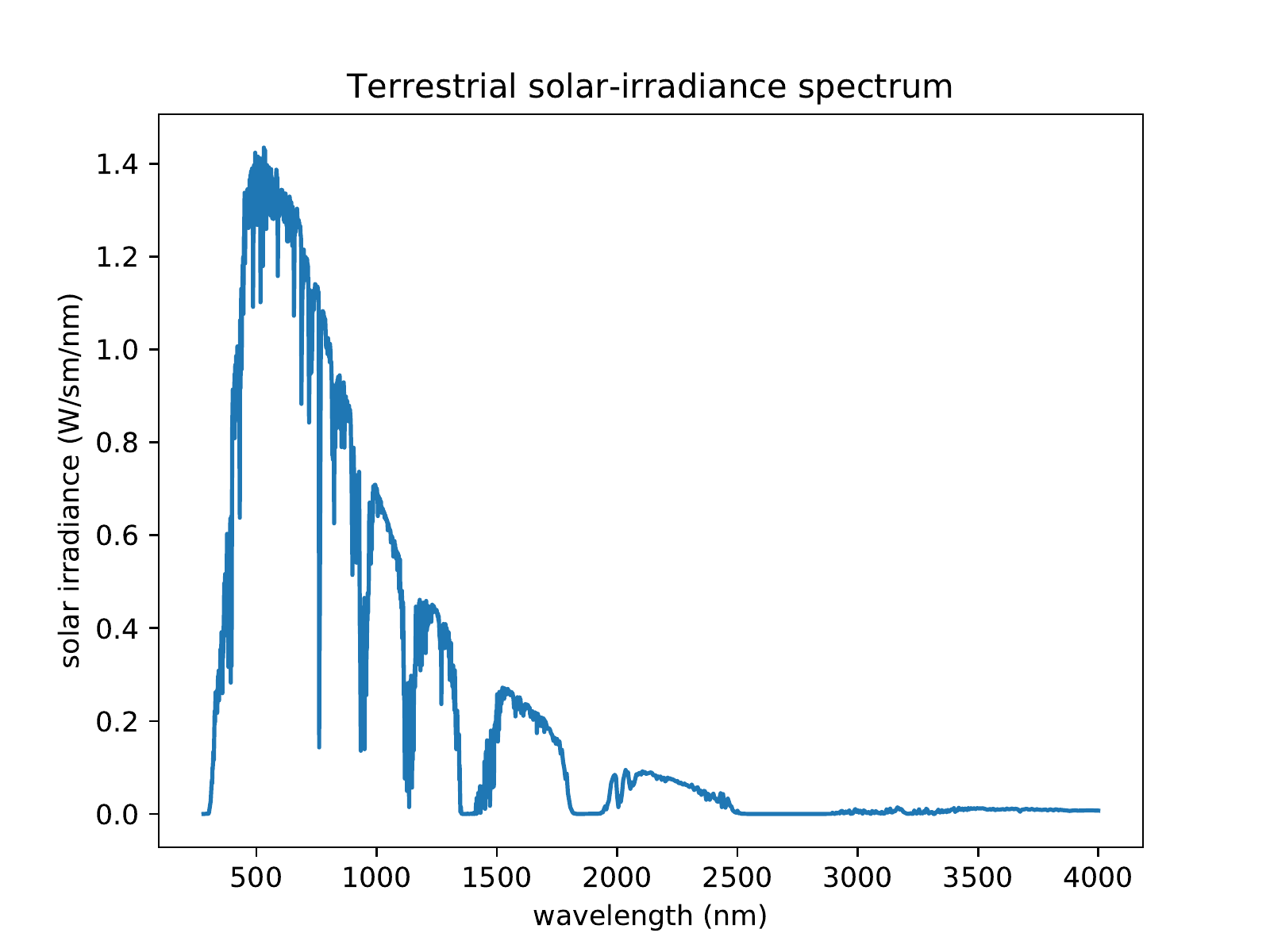}
\par\end{centering}
\caption{ASTM AM1.5 solar-irradiance spectrum~\cite{ASTM05}, given by 2002
tabulated data points for wavelengths $\lambda$ from 280--4000nm.
The spectrum is extremely oscillatory, with many sharp dips corresponding
to absorption by various molecules.\label{fig:irradiance}}

\end{figure}

For many problems involving solar energy, such as simulating the efficiency
of solar cells, it is necessary to compute integrals over wavelength~$\lambda$
of the form
\begin{equation}
\int_{a}^{b}S(\lambda)f(\lambda)d\lambda,\label{eq:Sintegral}
\end{equation}
where $S(\lambda)\ge0$ is the \emph{solar irradiance} spectrum (the
intensity of sunlight on the surface of the Earth), $(a,b)$ is some
bandwidth of interest, and $f(\lambda)$ is an application-specific
integrand. For example, in solar photovoltaic devices, $f(\lambda)$
is computed from the optical absorption spectrum of a device in order
to quantify the photonic efficiency~\cite{Bermel2007}. The basic
challenge in computing eq.~(\ref{eq:Sintegral}) is that $S(\lambda)$
is \emph{extremely oscillatory }and \emph{non-smooth }tabulated data
(given by the ASTM standard~\cite{ASTM05}, shown in Fig.~\ref{fig:irradiance}),
while $f(\lambda)$ is typically smooth but computationally expensive
(e.g. requiring a solution of Maxwell's equations for each $\lambda$~\cite{Bermel2007}).
 Worse, for device design and optimization the integral (\ref{eq:Sintegral})
must be computed many times for different $f(\lambda)$~\cite{Bermel2007,ShengJo11,Oskooi2012,Lin2013,Wiesendanger2013}.
In these notes, we explain how eq.~(\ref{eq:Sintegral}) can be accurately
approximated with only a small number of $f(\lambda)$ evaluations
by a \emph{Gaussian quadrature} rule~\cite{BrassPetras11}
\begin{equation}
\int_{a}^{b}S(\lambda)f(\lambda)d\lambda\approx\sum_{i=1}^{N}w_{i}f(\lambda_{i}),\label{eq:Squadrature}
\end{equation}
where $\lambda_{i}$ are the quadrature points (nodes), $w_{i}$ are
the quadrature weights, and $N$ is the \emph{order} of the rule.
The key fact is that we can construct a Gaussian quadrature rule in
which the effect of the complicated irradiance function $S(\lambda)$
is \emph{precomputed} in the $(w_{i},\lambda_{i})$ data, so that
any polynomial $f(\lambda)$ of degree up to $2N-1$ is integrated
\emph{exactly} by eq.~(\ref{eq:Squadrature}) and smooth $f(\lambda)$
are integrated with an error that decreases \emph{exponentially} with
$N$~\cite{BrassPetras11,Trefethen2008}. 

In the following notes, we first briefly review the numerical construction
of Gaussian quadrature rules (Sec.~\ref{sec:Constructing-Gaussian-quadrature}),
then present the resulting solar-weighted quadrature rules and demonstrate
their accuracy on some examples (Sec.~\ref{sec:Solar-integration-rules}),
and finally comment on various possible extensions and variations
(Sec.~\ref{sec:Extensions-and-variations}). We provide a selection
of precomputed rules~(\ref{eq:Squadrature}) as well as a program
in the Julia language~\cite{bezanson2017julia} (in the form of a
Jupyter notebook~\cite{Kluyver:2016aa} tutorial) to compute a rule~(\ref{eq:Squadrature})
for any desired $N$ and bandwidth $(a,b)$~\cite{solarnotebook}.

\section{Constructing Gaussian-quadrature rules\label{sec:Constructing-Gaussian-quadrature}}

The construction of Gaussian quadrature rules, by which polynomial
$f(\lambda)$ are integrated exactly up to a given degree, is closely
tied to the theory of orthogonal polynomials~\cite{BrassPetras11}.
For a given weight function $S(x)\ge0$ and interval $(a,b)$, we
define an inner product 
\begin{equation}
\langle p_{1},p_{2}\rangle=\int_{a}^{b}S(x)p_{1}(x)p_{2}(x)dx\label{eq:inner}
\end{equation}
of any polynomial functions $p_{1,2}(x)$. Then, by a Gram--Schmidt
orthogonalization process applied to $\{1,x,x^{2},\ldots\}$, one
can obtain a sequence of orthonormal polynomials $\{q_{0}(x),q_{1}(x),\ldots\}$
with respect to this inner product, where $q_{n}(x)$ has degree~$n$.
Remarkably, it turns out that the quadrature points $\lambda_{i}$
in eq.~(\ref{eq:Squadrature}) are exactly the roots of $q_{N}(x)$,
and the weights can also be computed from these polynomials~\cite{BrassPetras11}.

Numerically, there are a variety of schemes for computing the quadrature
points $\lambda_{i}$ and weights $w_{i}$. One of the simplest accurate
methods is the Golub--Welsch algorithm~\cite{Golub1969}, which
constructs a tridiagonal ``Jacobi'' matrix corresponding to a three-term
recurrence for the $q_{k}$ polynomials, and then obtains the nodes
$\lambda_{i}$ and weights $w_{i}$ from the eigenvalues and eigenvectors
of this matrix, respectively. It turns out that this procedure is
equivalent to a Lanczos iteration for the symmetric linear operator
corresponding to multiplication by $x$~\cite{Trefethen1997}. The
only input that is required is a way to compute the functional
\begin{equation}
S\{p\}=\int_{a}^{b}S(x)p(x)dx,\label{eq:Spoly}
\end{equation}
which is just the integral~(\ref{eq:Sintegral}) computed for any
given polynomial $p(x)$ (up to degree $2N$). For the most famous
applications of Gaussian quadrature to simple weight functions like
$1$ on $(-1,1)$ (Gauss--Legendre quadrature), $e^{-x^{2}}$ on
$(-\infty,+\infty)$ (Gauss--Hermite quadrature), or $e^{-x}$ on
$(0,\infty)$ (Gauss--Laguerre quadrature), these integrals and the
resulting three-term recurrences are known analytically. For an arbitrary
weight function $S(x)$ like the solar-irradiance spectrum, we must
perform integral~(\ref{eq:Spoly}) numerically (for $2N$ polynomials
in total).

This process may at first seem fruitless: isn't eq~(\ref{eq:Spoly})
equivalent to our original problem? There are two critical differences,
however. First, polynomials $p(x)$ are very \emph{cheap to evaluate},
so we can afford to use brute-force numerical integration methods
that require $\sim10^{4}$ function evaluations to evaluate~(\ref{eq:Spoly}),
unlike solar-cell applications of eq.~(\ref{eq:Sintegral}) where
the integrand $f(\lambda)$ is very expensive. Second, we only need
to perform these polynomial integrations \emph{once} for a given $N$:
we can then \emph{re-use} the resulting quadrature rule $(\lambda_{i},w_{i})$
over and over for many different $f(\lambda)$.

\begin{figure}
\begin{centering}
\includegraphics[width=0.7\columnwidth]{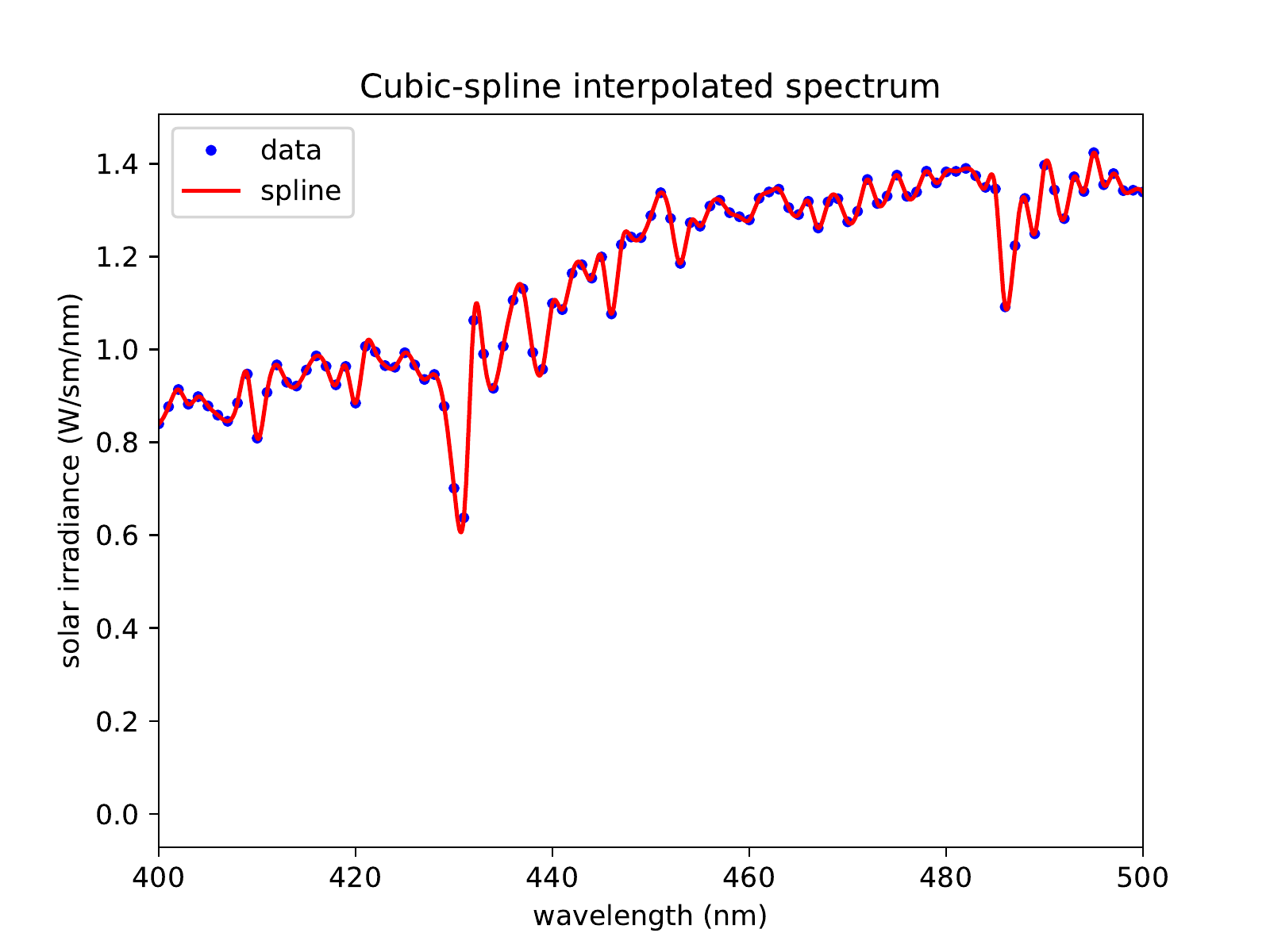}
\par\end{centering}
\caption{Cubic-spline interpolation~\cite{Dierckx93} of the ASTM AM1.5 solar-irradiance
spectrum~\cite{ASTM05}, shown for a portion of the spectrum. To
ensure positivity of the interpolant for all $\lambda$, we actually
fit a spline $C(\lambda)$ to the square root of the irradiance, and
approximate $S(\lambda)$ by $C(\lambda)^{2}$ for all $\lambda\in(280,4000)$nm.\label{fig:spline}}
\end{figure}

In the specific case of our solar-irradiance spectrum $S(\lambda)$,
the data are supplied in the form of 2002 data points from 280nm to
4000nm by the ASTM standard~\cite{ASTM05}. To define $S(\lambda)$
continuously for all $\lambda$ in this range, we perform a cubic-spline
interpolation of the data using a library by Dierckx~\emph{et al.~}\cite{Dierckx93}.
More precisely, since we require an interpolant $S(\lambda)\ge0$
everywhere in order for eq.~(\ref{eq:inner}) to define a proper
inner product, we compute a cubic-spline interpolation $C(\lambda)$
of $\sqrt{S(\lambda)}$ at the ASTM data points, and then approximate
$S(\lambda)$ by $C(\lambda)^{2}$ everywhere. This spline interpolant
is shown for a portion of the spectrum from 400--500nm in Fig.~\ref{fig:spline}.
 A cubic spline has the property that not only does it pass through
all of the data points, but its first and second derivatives are also
continuous; the data points where the spline switches from one cubic
curve to another (where the third derivative is discontinuous) are
known as the \emph{knots} of the spline~\cite{Dierckx93}. In order
to perform the $S\{p\}$ integrals~(\ref{eq:Spoly}) of polynomials
$p(\lambda)$ against $C(\lambda)^{2}$, we simply break the integral
into $\approx2000$ segments at the knots (to avoid integrating through
discontinuities) and apply a globally adaptive~\cite{Malcolm1975}
Gauss--Kronrod~\cite{BrassPetras11} quadrature scheme implemented
in Julia~\cite{quadgk} to high accuracy ($>$ 9 digits).

Given an arbitrary weight function (and optionally a specialized routine
to compute the $S\{p\}$ integrals), the QuadGK package in Julia~\cite{quadgk}
provides a subroutine to execute the Golub--Welsch algorithm and
compute the quadrature points and weights. The only addition to the
textbook description of the algorithm~\cite{Trefethen1997}, besides
the numerical computation of $S\{p\}$, is that the integration domain
$(a,b)$ is internally rescaled to $(-1,1)$ so that the polynomials
$p(x)$ can be represented in the basis of the Chebyshev polynomials
$\{T_{0},T_{1},\ldots\}$~\cite{Trefethen12}. Given the coefficients
$c_{k}$ of a Chebyshev series $p(x)=\sum_{j}a_{k}T_{k}(x)$, the
polynomial $p(x)$ is evaluated by a Clenshaw recurrence~\cite{Trefethen12},
and another recurrence is used to compute the coefficients of the
$xp(x)$ polynomial required for the Golub--Welsch algorithm~\cite{Trefethen1997}.
This approach\footnote{A similar Chebyshev-based Lanczos recurrence was implemented by S.~Olver
in 2014 for the ApproxFun Julia package~\cite{Olver2014} following
a suggestion by A.~Edelman.} to representing and evaluating polynomials avoids the numerical problems
(ill conditioning) that arise in the familiar monomial basis $\{1,x,x^{2},\ldots\}$---Chebyshev
polynomials $T_{k}(x)=\cos(k\cos^{-1}x)$ are well-behaved numerically
since they are ``cosines in disguise'' \cite{Trefethen12}. (A related,
but more complicated, algorithm in which one first computes the ``moments''
$S\{T_{k}(x)\}$, is reviewed by Gautschi~\cite{Gautschi1994}, and
various other methods have also been proposed~\cite{Gautschi1994,Mantica1996,Gander2001,Fernandes2006}.)

\section{Solar integration rules and convergence\label{sec:Solar-integration-rules}}

\begin{figure}
\begin{centering}
\includegraphics[width=0.7\columnwidth]{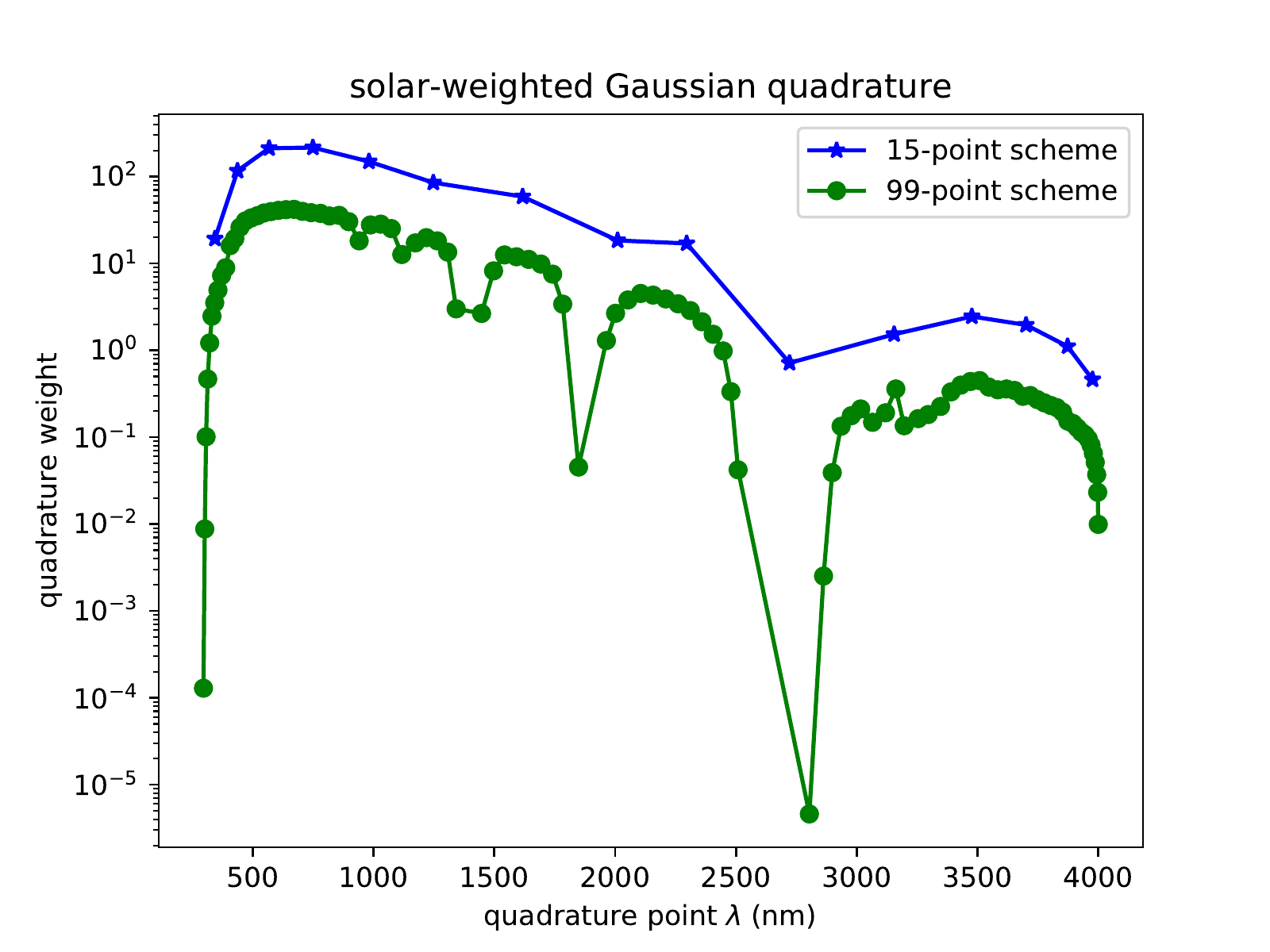}
\par\end{centering}
\caption{\label{fig:quad}Solar-weighted Gaussian quadrature rules (ASTM AM1.5
spectrum~\cite{ASTM05}) in a bandwidth of 280--4000nm, for orders
$N=15$ and $N=99$.}

\end{figure}
Given the computational methods in Sec.~\ref{sec:Constructing-Gaussian-quadrature},
we can then construct solar-weighted Gaussian quadrature rules for
any desired order $N$ and any given bandwidth $\lambda\in(a,b)\subseteq(280,4000)\text{nm}$
using the code provided in the attached Julia notebook~\cite{solarnotebook}.
The resulting weights $w_{i}$ are plotted versus the corresponding
quadrature points $\lambda_{i}$ for $N=15$ and $N=99$. Not surprisingly,
they roughly follow the peaks and dips of the solar spectrum $S(\lambda)$
from Fig.~\ref{fig:irradiance}.

\begin{figure}

\begin{centering}
\includegraphics[width=0.7\columnwidth]{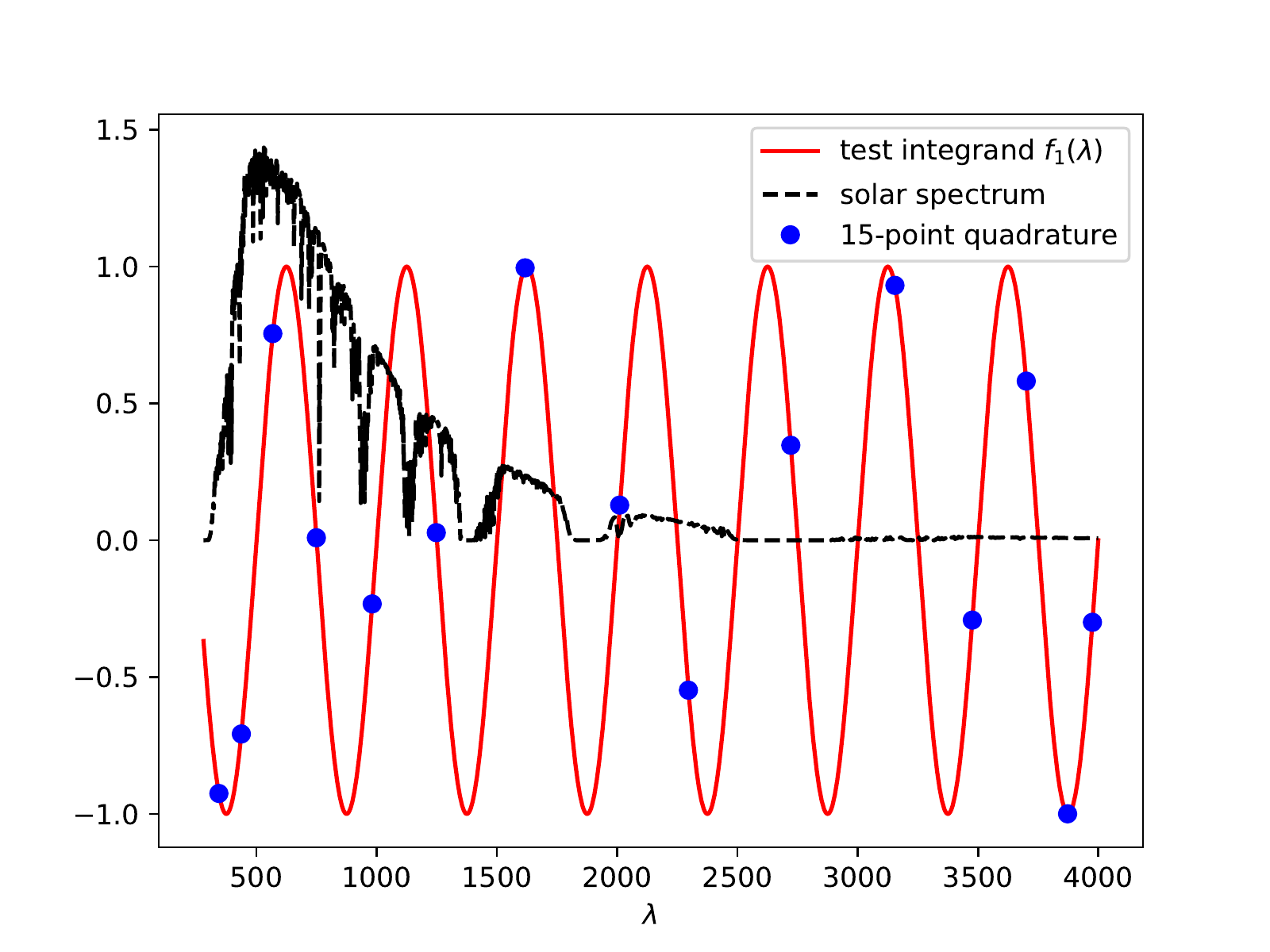}
\par\end{centering}
\caption{Test integrand $f_{1}(\lambda)=\sin(2\pi\lambda/500)$, superimposed
on the solar-irradiance spectrum $S(\lambda)$. The blue dots show
$f_{1}(\lambda_{i})$ at the quadrature points $\lambda_{i}$ of the
$N=15$ point solar-weighted Gaussian quadrature rule from Fig.~\ref{fig:quad}.\label{fig:test1}}

\end{figure}
To validate these quadrature schemes and evaluate their accuracy,
let us consider two test integrands, $f_{1}(\lambda)=\sin(2\pi\lambda/500)$
and $f_{2}(\lambda)=\sin(2\pi\lambda/50)$. The first integrand $f_{1}$
is shown in Fig.~\ref{fig:test1}: it oscillates 8 times over the
solar bandwidth, and seems to be barely sampled adequately by the
$N=15$ point quadrature rule. To compute the ``exact'' integrals
of $f_{1}$ and $f_{2}$, we apply the same method that we used for
numerical polynomial integration $S\{p\}$ in Sec.~\ref{sec:Constructing-Gaussian-quadrature}:
we partition the integral into $\approx2000$ segments at the knots
of the cubic-spline interpolant for $S(\lambda)$ and then apply an
adaptive Gauss--Kronrod integration scheme to nearly machine precision
($>12$ significant digits). Comparing this brute-force result to
that of the quadrature rules~(\ref{eq:Squadrature}) applied to $f_{1}(\lambda)$
from Fig.~\ref{fig:test1}, we find that the $N=15$ quadrature rule
yields the correct result with an error of only $0.7\%$, while the
$N=99$ quadrature rule yields the correct result to at least 13 significant
digits.

\begin{figure}
\begin{centering}
\includegraphics[width=0.7\columnwidth]{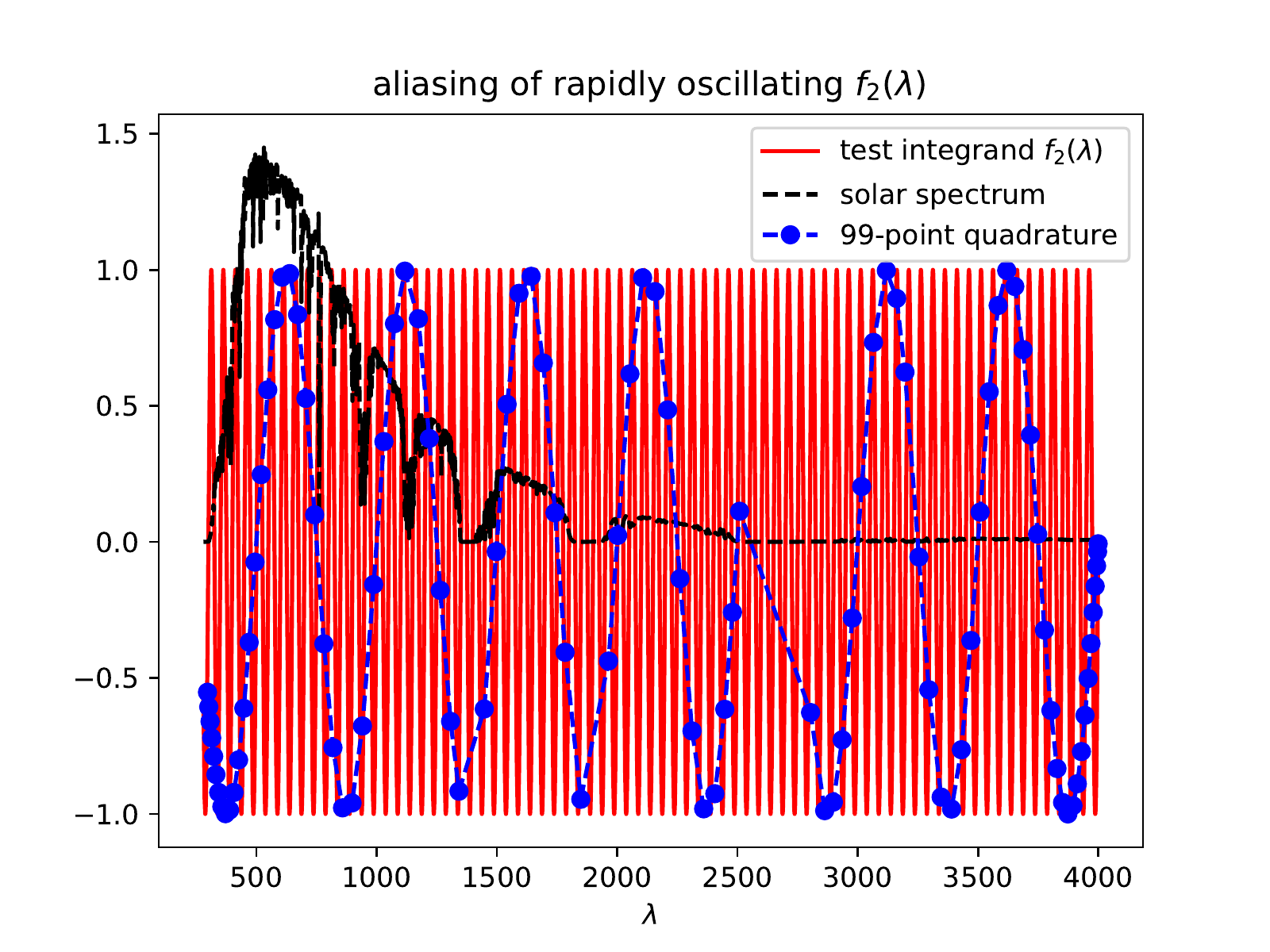}
\par\end{centering}
\caption{Test integrand $f_{2}(\lambda)=\sin(2\pi\lambda/50)$, superimposed
on the solar-irradiance spectrum $S(\lambda)$. The blue dots show
$f_{2}(\lambda_{i})$ at the quadrature points $\lambda_{i}$ of the
$N=99$ point solar-weighted Gaussian quadrature rule from Fig.~\ref{fig:quad}.
This 99-point sampling is too coarse to capture the rapid oscillation
of $f_{2}$, and results in an aliasing phenomenon~\cite{Trefethen2008}
in which the sampled function appears to oscillate at a lower frequency.\label{fig:test2}}
\end{figure}
For the test function $f_{2}(\lambda)$, which oscillates $10\times$
more rapidly than $f_{1}$, even the $N=99$ quadrature rule is not
enough to obtain an accurate result: the sampling of $f_{2}$ is simply
not fine enough to detect its rapid variations. This can be seen in
Fig.~\ref{fig:test2}, which depicts the function $f_{2}$ sampled
at the quadrature points of the $N=99$ rule. The sampling clearly
exhibits the phenomenon of \emph{aliasing}, which is central to understanding
the convergence rates of these quadrature schemes~\cite{Trefethen2008}:
the sampled function $f_{2}(\lambda_{i})$ appears to oscillate at
a much lower frequency. Since aliasing means that the quadrture rule
cannot ``perceive'' the actual rapid oscillation of $f_{2}$, the
quadrature formula~(\ref{eq:Squadrature}) unsurprisingly yields
a completely incorrect result ($>100$\% error). A larger $N$ must
be employed to integrate such a rapidly oscillating $f(\lambda)$.
In the case of this $f_{2}(\lambda)$, an $N=140$ solar-weighted
Gaussian quadrature rule is sufficient to compute $\int S(\lambda)f_{2}(\lambda)$
to 9 significant digits. The key point is that the required order
of the quadrature rule is completely determined by the features of
the application-supplied $f(\lambda)$, and is \emph{independent}
of the preprocessed complexity of the solar irradiance $S(\lambda)$.

\begin{figure}
\begin{centering}
\includegraphics[width=0.7\columnwidth]{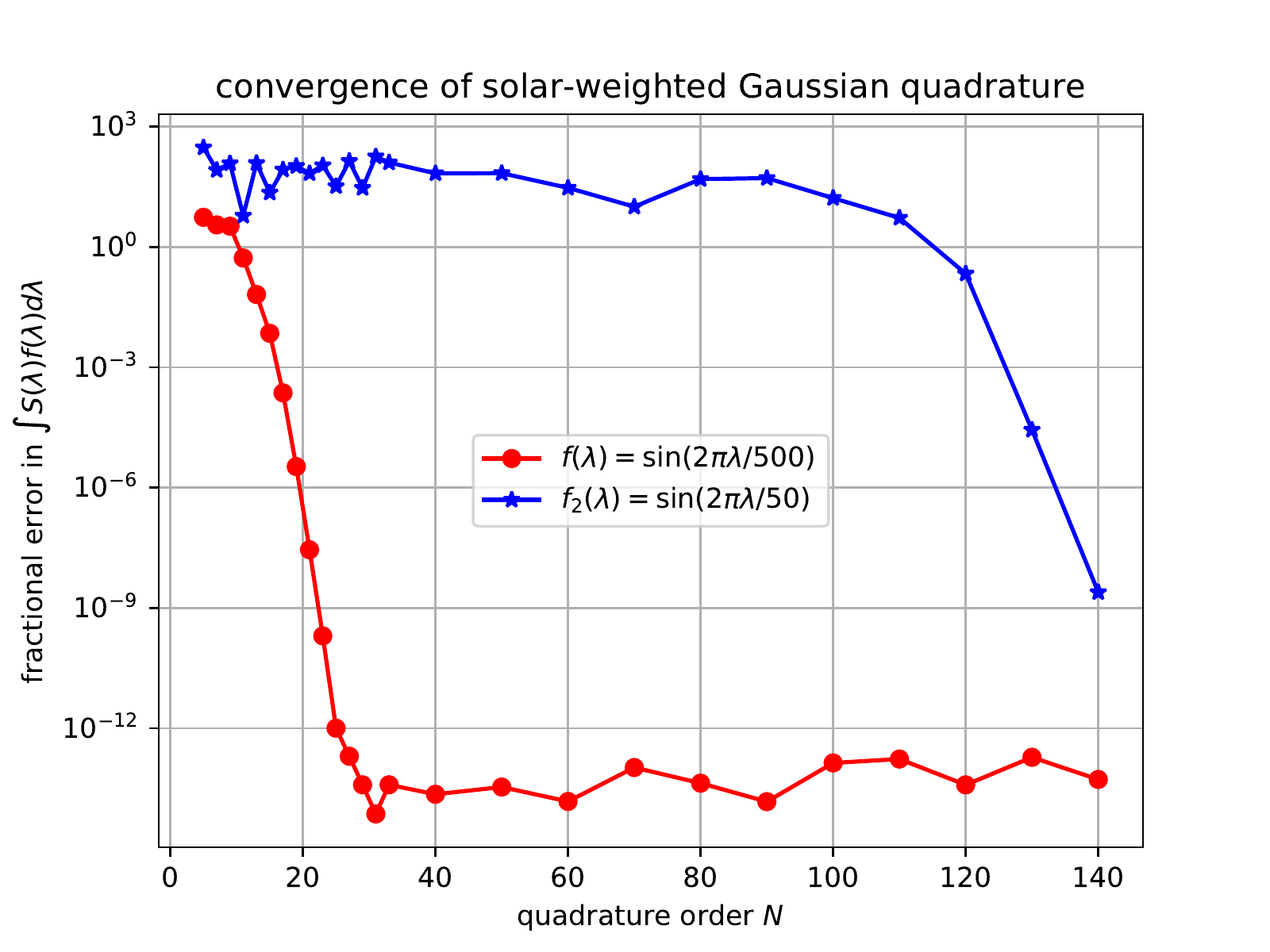}
\par\end{centering}
\caption{Convergence of solar-weighted Gaussian quadrature: fraction error
in the integral $\int S(\lambda)f(\lambda)d\lambda$ versus the number
$N$ of quadrature points for the two test integrands $f_{1}(\lambda)=\sin(2\pi\lambda/500)$
and $f_{2}(\lambda)=\sin(2\pi\lambda/50)$. In both cases, once $N$
becomes sufficiently large, the error falls exponentially with $N$.
(However, in the case of the rapidly oscillating $f_{2}$ from Fig.~\ref{fig:test2},
$N$ must first become $>100$ in order to adequately sample $f_{2}$'s
oscillations.) The error in the $f_{1}$ does not drop below a ``noise
floor'' of $\approx10^{-13}$ due to the limitations of double-precision
floating-point arithmetic~\cite{Trefethen1997}.\label{fig:convergence}}
\end{figure}
In general, when Gaussian quadrature is applied to any smooth (infinitely
differentiable) function $f(\lambda)$, such as an absorption spectrum,
the error converges to zero faster than any power of $1/N$; typically
(for analytic functions), the error converges \emph{exponentially}
with $N$~\cite{Trefethen2008}. This rapid convergence\footnote{In fact, since $f_{1}$ and $f_{2}$ are ``entire'' functions analytic
everywhere in the complex plane, the convergence should asymptotically
be \emph{faster} than exponential in $N$. In practice, this is hard
to detect because the accuracy is quickly constrained by machine precision.} is demonstrated in Fig.~\ref{fig:convergence} for both $f_{1}$
and $f_{2}$. The accuracy of the $f_{1}$ integrand, in fact, nearly
hits the limits of machine precision around $N\approx30$ (and might
be able to go slightly lower than $10^{-13}$ if we reduced the tolerances
in the brute-force integration used in this section and in Sec.~\ref{sec:Constructing-Gaussian-quadrature}).
The $f_{2}$ integration accuracy also exhibits exponential convergence,
but the error does not begin its asymptotic decline until $N\gtrsim100$,
at which point the quadrature rule begins to adequately sample $f_{2}$'s
rapid oscillation without destructive aliasing.

\section{Extensions and variations\label{sec:Extensions-and-variations}}

\begin{figure}
\begin{centering}
\includegraphics[width=0.7\columnwidth]{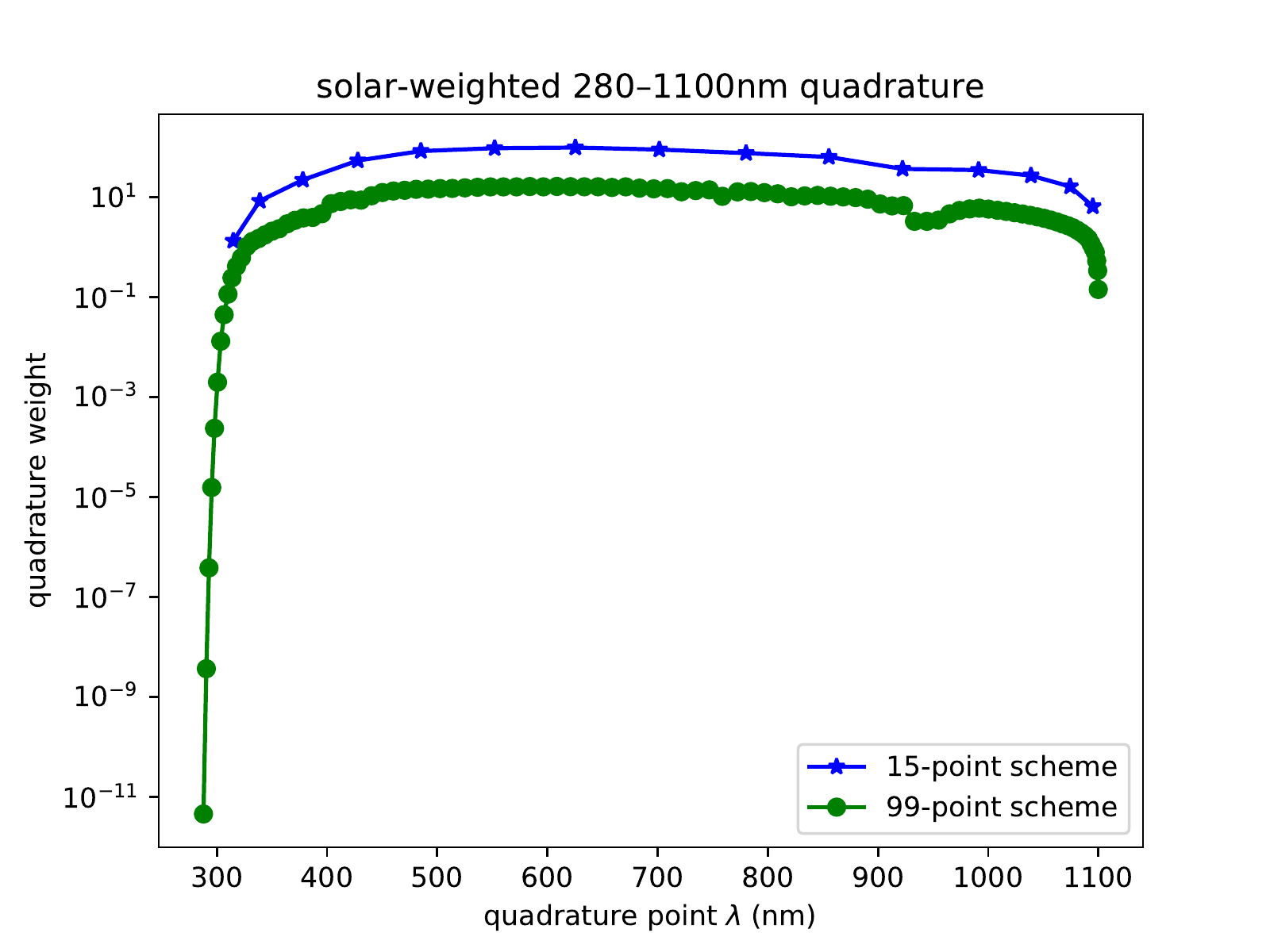}
\par\end{centering}
\caption{\label{fig:quad2}Solar-weighted Gaussian quadrature rules (ASTM AM1.5
spectrum~\cite{ASTM05}) in a bandwidth of 280--1100nm, for orders
$N=15$ and $N=99$.}
\end{figure}
To apply the solar-weighted quadrature rule to specific applications,
one may wish to consider some variations. For example, with silicon
solar cells one rarely considers wavelengths beyond 1100nm, since
longer wavelengths fall short of the bandgap of silicon and cannot
easily generate electron--hole pairs~\cite{Bermel2007,Oskooi2012}.
So, for these devices one might want a quadrature rule for $\int_{280\text{nm}}^{1100\text{nm}}S(\lambda)f(\lambda)d\lambda$.
Two such quadrature rules, for $N=15$ and $N=99$ points, are shown
in Fig.~\ref{fig:quad2}. It is advantageous to specialize the quadrature
rule to the bandwidth of interest because one can then resolve much
finer features in $f(\lambda)$ for the same number of points, and
optimized solar cells typically exhibit many sharp resonant peaks~\cite{Bermel2007,ShengJo11,Oskooi2012}.
As another example, the figure of merit for a photovoltaic cell is
typically proportional to $\int S(\lambda)\lambda A(\lambda)d\lambda$
where $A(\lambda)$ is the device's absorption spectrum~\cite{Bermel2007},
so one could also include the factor of $\lambda$ in the precomputed
weight when constructing the Gaussian-quadrature schemes. It might
also be useful to rescale the weight by the absorption coefficient
of the photovoltaic material, to further skew the quadrature points
towards the wavelengths where more absorption occurs.

Naturally, the same techniques could also be applied to completely
different weight functions. Although weighted Gaussian quadrature
schemes have been well known for decades in the numerical-integration
literature, they seem to have been rarely applied for the complicated
weighting functions that often arise in optical physics. For example,
the CIE color-projection functions (used to model the perception of
a human eye)~\cite{Schanda07} are described by tabulated data (though
there are also analytical fits~\cite{Wyman2013xyz}), leading to
integrals that qualitatively resemble Gauss--Hermite quadrature but
require a quantitatively distinct quadrature scheme for high accuracy.
Thermal emission from heated bodies is another good example: the emitted
power is an integral of a surface's absorption spectrum (found by
solving Maxwell's equations) weighted by a black-body ``Planck''
spectrum $\sim\omega^{3}/(e^{\hbar\omega/kT}-1)$~\cite{Cornelius1999},
suggesting a need for a Gauss--Laguerre-like Planck-weighted quadrature
scheme.

\section*{Acknowledgements}

This work was supported in part by the U. S. Army Research Office
through the Institute for Soldier Nanotechnologies under grant W911NF-13-D-0001.

\bibliographystyle{IEEEtran}
\bibliography{solar}

\end{document}